\documentclass{amsart}

\usepackage{palatino}
\usepackage{amssymb,mathrsfs,latexsym,amssymb}

\newcommand{\SP}{\operatorname{Sp}}
\newcommand{\GL}{\operatorname{GL}}

\newcommand{\SL}{\operatorname{SL}}

\newcommand{\End}{\operatorname{End}}
\newcommand{\tensor}{\otimes}

\newcommand{\stab}{\operatorname{Stab}}
\newcommand{\ch}{\operatorname{ch}}

\newcommand{\ext}{\operatorname{Ext}}

\newcommand{\res}{\operatorname{res}}

\newcommand{\XX}{\mathscr{X}}
\newcommand{\YY}{\mathscr{Y}}
\newcommand{\ZZ}{\mathscr{Z}}

\newcommand{\pow}[1]{[ \hspace{-1.6pt} [ {#1} ] \hspace{-1.6pt} ]}

\DeclareMathAlphabet{\mathrc}{U}{eur}{m}{n}

\newcommand\C{\mathbf{C}}    %        complex field
\newcommand\Z{\mathbf{Z}}    %        integers
\newcommand\F{\mathbf{F}}    %        finite fields
\newcommand\G{\mathbf{G}}    %        mult or additive group
\newcommand\Lie{\operatorname{Lie}}
\newcommand\Ad{\operatorname{Ad}}
\newcommand\ad{\operatorname{ad}}

\newcommand\lie[1]{\mathfrak{#1}}
\newcommand\glie{\lie{g}}
\newcommand\olie{\lie{o}}

\newcommand\blie{\lie{b}}
\newcommand\slie{\lie{s}}
\newcommand\plie{\lie{p}}
\newcommand\gl{\lie{gl}}

\newcommand{\hlie}{\lie{h}}
\newcommand{\ulie}{\lie{u}}

\newcommand{\sllie}{\lie{sl}}
\newcommand{\clie}{\lie{c}}

\newcommand{\nlie}{\lie{n}}

\newcommand{\bwedge}{\textstyle{\bigwedge}}

\newcommand{\congruent}{\equiv} 
\newcommand{\iso}{\simeq}

\renewcommand{\mathbb}{\mathbf}

\newcommand{\e}{\varepsilon}

\newcommand{\CC}{\mathscr{C}}

\newcommand{\NN}{\mathcal{N}}
\newcommand{\UU}{\mathcal{U}}

\newcommand{\witt}{\mathfrak{w}}

\newcommand{\gr}{\operatorname{gr}}

\newcommand{\formal}{\mathscr{F}}
\newcommand{\Sym}{\operatorname{Sym}}

\newcommand{\tensorpoly}{\textstyle{\bigotimes}}

\newtheorem{theorem}{Theorem}
\newtheorem{prop}[theorem]{Proposition}
\newtheorem{cor}[theorem]{Corollary}
\newtheorem{lem}[theorem]{Lemma}
\theoremstyle{remark}
\newtheorem{rem}[theorem]{Remark}

\newtheorem{example}[theorem]{Example}

\numberwithin{equation}{section}

\begin{document}

\begin{abstract}
  Let $G$ be a quasisimple algebraic group over an algebraically
  closed field of characteristic $p>0$. We suppose that $p$ is
  \emph{very good} for $G$; since $p$ is good, there is a bijection between the
  nilpotent orbits in the Lie algebra and the unipotent classes in
  $G$. If the nilpotent $X \in \Lie(G)$ and the unipotent $u \in G$
  correspond under this bijection, and if $u$ has order $p$, we show
  that the partitions of $\ad(X)$ and $\Ad(u)$ are the same.
  When $G$ is classical or of type $G_2$, we prove this result with no
  assumption on the order of $u$.
  
  In the cases where $u$ has order $p$, the result is achieved through
  an application of results of Seitz concerning good $A_1$ subgroups
  of $G$. For classical groups, the techniques are more
  elementary, and they lead also to a new proof
  of the following result of Fossum: the structure constants of the
  representation ring of a 1-dimensional formal group law $\formal$
  are independent of $\formal$.
\end{abstract}

\bibliographystyle{hamsalpha}
\title{Adjoint Jordan blocks}
\author{George McNinch}
\thanks{This work was supported in part by a grant from the National
  Science Foundation.}

\address{Department of Mathematics \\
         Room 255 Hurley Building \\
         University of Notre Dame \\
         Notre Dame, Indiana 46556-5683 \\
         USA}
\email{mcninch.1@nd.edu}
\date{June 28, 2002}

\maketitle

\section{Introduction}

Let $G$ be a quasisimple algebraic group over an algebraically closed
field $k$ of characteristic $p>0$. In case $p$ is \emph{good} for
$G$, there is a bijection between the nilpotent orbits in the Lie
algebra, and the unipotent classes in $G$.
  
Suppose the class of the unipotent $u$ and the nilpotent $X$
correspond under this bijection. Assume further that $u$ has order $p$
(or equivalently, that $X$ is $p$-nilpotent, i.e. $X^{[p]} = 0$). We
prove in Theorem \ref{theorem:main-theorem-p-nilpotent} that if $p$ is
\emph{very good}, the partition of $\Ad(u)$ is the same as the
partition of $\ad(X)$ (where $\Ad$ and $\ad$ denote the respective
adjoint actions of $G$ and $\glie = \Lie(G)$ on $\glie$).  This is
carried out in \S \ref{sec:good_a1s} using results of Seitz
\cite{seitz-unipotent} on good $A_1$-type subgroups of $G$.

When $G$ is a classical group, or has type $G_2$, we show in Theorems
\ref{theorem:classical-adjoint-partitions} and
\ref{theorem:g2-adjoint-partitions} that the partitions of $\Ad(u)$
and $\ad(X)$ coincide with no assumption on the order of $u$.  The
techniques for classical groups involve comparison of ``group like''
and ``Lie algebra like'' tensor products; see \S \ref{sec:similar}.
One may consider more generally the tensor product determined by any 1
dimensional formal group law.  In this context, our techniques provide
a new proof of a theorem of Fossum: the representation ring of a
formal group law $\formal$ is independent of $\formal$; see Corollary
\ref{cor:rep-ring}.

Lawther \cite{Lawther:Jblock} has found (with the aid of a computer)
the adjoint partitions of the unipotent classes in exceptional simple
algebraic groups. Combined with his work, Theorem
\ref{theorem:main-theorem-p-nilpotent} computes the adjoint partitions
for those nilpotent classes of exceptional groups in good characteristic
which are $p$-nilpotent.

In a final section \S \ref{sec:char0}, we consider adjoint partitions
in characteristic 0. Let $G$ be a semisimple group of rank $r$ in
characteristic 0.  We show that for certain nilpotent elements $X$,
the eigenvalues (on a Cartan subalgebra) of a corresponding Weyl group
element account for $r$ parts of the partition of $\ad(X)$.

At the risk of being pedantic, we recall: The partition of a nilpotent
linear endomorphism $X$ of a vector space $V$ over a field is the
ordered sequence of Jordan block sizes of $X$.  The partition of a
unipotent linear automorphism $u$ of $V$ is the partition of the
nilpotent map $u-1$.

\section{Simple groups}

Let $G$ be a quasisimple algebraic group over $k$ with Lie algebra
$\glie$; thus the root system $R$ of the semisimple group $G$ is
irreducible.

\subsection{The bad, the good, and the very good}

  The characteristic $p$ of $k$ is said to be bad for $R$
in the following circumstances: $p=2$ is bad whenever $R \not = A_r$,
$p=3$ is bad if $R = G_2,F_4,E_r$, and $p=5$ is bad if $R=E_8$.
Otherwise, $p$ is good. 

Moreover, $p$ is very good if $r \not \congruent -1 \pmod p$ in case
$R=A_r$.

\subsection{Classical groups}
\label{sub:classical}
Let $V$ be a finite dimensional vectorspace over $k$, suppose that
$\varphi$ is a bilinear form on $V$, let $\Omega$ be the full group of
isometries of $V$ with respect to $\varphi$, and let $\olie =
\Lie(\Omega)$.  We will say that $\Omega$ is a \emph{classical group}
in the following situations:
\begin{enumerate}
\item[CG1.] $\varphi = 0$, so that $\Omega = \GL(V)$ and $\olie = \gl(V)$.
\item[CG2.] $\varphi$ is nondegenerate and alternating, so that $\Omega =
  \SP(V,\varphi)$ is a symplectic group, and $\olie =
  \lie{sp}(V,\varphi)$.
\item[CG3.] $\varphi$ is nondegenerate and symmetric, so that $\Omega =
  O(V,\varphi)$ is an orthogonal group, and $\olie =
  \lie{so}(V,\varphi)$.
\end{enumerate}

In each case, $\Omega$ is a reductive group (though it is not
connected in case CG3, and it is not semisimple in case CG1).  The
prime 2 is bad in cases CG2 and CG3; all other primes are good for
classical groups.  For convenience, we write $(\nu,V)$ for the natural
representation of $\Omega$.

We record the following well-known characterization of the adjoint
representation of a classical group $\Omega$.
\begin{lem}
  \label{lem:classical-adjoint-modules}
  There is an isomorphism of
  $\Omega$-representations
  \begin{equation*}
    (\Ad,\olie) \iso \left \{ 
      \begin{matrix}
        (\nu \tensor \nu^\vee,V \tensor V^\vee) &  \text{in case CG1,}\\
        (\Sym^2\nu,\Sym^2 V) &  \text{in case CG2,} \\
        (\bigwedge^2 \nu,\bwedge^2 V) &  \text{in case CG3.} \\
      \end{matrix}
    \right .
  \end{equation*} 
\end{lem}
Here, $V^\vee$ denotes the dual vector space (and $\nu^\vee$ the
contragredient representation).  
We also record:
\begin{lem}
  \label{lem:tensor-square}
  Let $W$ be a $k$-vector space. If the characteristic of $k$ is not
  2, there is an isomorphism of $\GL(W)$-representations
  \begin{equation*}
    W \tensor W \iso \Sym^2 W \oplus \bwedge^2 W.
  \end{equation*}
\end{lem}
 
\subsection{Springer's isomorphism}

Let $G$ be connected, reductive with an irreducible root system, let $\NN \subset
\glie$ be the nilpotent variety , and let $\UU \subset G$ be the
unipotent variety of $G$.  We recall
the following result (which depends on Bardsley and Richardson's
construction of a ``Springer isomorphism'').
\begin{prop} 
  \label{prop:equivariant-iso}
  If the characteristic of $k$ is very good for $G$, there is a
  $G$-equivariant isomorphism $\theta:\NN \to \UU$ with
  the property  $\theta(X^{[p]}) = \theta(X)^p$ for each
  $X \in \NN$.
\end{prop}

\begin{proof}
  \cite[Lemma 27, Theorem 35]{Mc:sub-principal}.
\end{proof}

\subsection{Distinguished nilpotents and associated co-characters}
\label{sub:distinguished}

In this paragraph, $G$ may be an arbitrary reductive group in good
characteristic. A nilpotent $X \in \glie$ is \emph{distinguished} if
the connected center of $G$ is a maximal torus of $C_G(X)$. The
Bala-Carter theorem, proved by Pommerening for reductive groups in
good characteristic, implies that any nilpotent $X \in \glie$ is
distinguished in the Lie algebra of some Levi subgroup $L$ of $G$.

Let $X \in \glie$ be nilpotent.  A co-character $\phi:k^\times \to G$
is said to be associated to $X$ if $\Ad \circ \phi(t)X = t^2X$ for all
$t\in k^\times$, and if the image of $\phi$ is contained in the
derived group $L'$ of a Levi subgroup $L$ such that $X$ is
distinguished in $\Lie(L)$.  The existence of a co-character
associated to a nilpotent $X \in \glie$ follows from Pommerening's
proof of the Bala-Carter theorem (more precisely: is a key step in his
proof of that theorem).  Any two co-characters associated with $X$ are
conjugate via $C_G^o(X)$; see e.g. \cite[Lemma
5.2]{jantzen:Nilpotent}.

The following useful observation is a consequence of definitions:
\begin{lem}
  \label{lem:sl2-disting}
  Let $\phi:\SL_{2/k} \to G$ be a homomorphism and suppose that $X=d\phi(
  \begin{pmatrix}
    0 & 1 \\
    0 & 0 
  \end{pmatrix})
  $ is a distinguished nilpotent element of $G$. Then the rule
  $\psi(t) = \phi(\text{diag}(t,t^{-1}))$ yields a co-character associated
  to $X$.
\end{lem}

\section{Adjoint Jordan blocks for elements of order $p$}

\subsection{Partitions and $\SL_2(k)$.}
\label{sec:good_a1s}

Let $S = \SL_2(k)$, and write $S(p) = \SL_2(\F_p)$.  The dominant
weights of a maximal torus of $S$ may be identified with $\Z_{\ge 0}$.
To each $d \in \Z_{\ge 0}$, there correspond various
$S$-representations: a simple module $L(d)$, a Weyl module $V(d)$ (of
dimension $d+1$), and an indecomposable tilting module $T(d)$; all
three have highest weight $d$. If $d < p$, then $L(d) = V(d) = T(d)$.
We refer to \cite{Donk3} for more details concerning tilting modules.
Note especially that each tilting $S$-module is a direct sum $\bigoplus_{d_i} T(d_i)$
for various $d_i \ge 0$; see \cite[Theorem 1.1]{Donk3}.

\begin{prop} Suppose that $p \not = 2$, and  $p \le d \le 2p-2$.
  \label{prop:tilting-partitions}
  \begin{enumerate}
  \item Each  unipotent element $1 \not = u \in S$ acts on $T(d)$
    with partition $(p,p)$.
  \item Each  nilpotent element $0 \not = X \in \Lie(S) = \lie{sl}_2(k)$
    acts on $T(d)$ with partition $(p,p)$.
  \end{enumerate}
\end{prop}

\begin{proof}
  By \cite[Lemma 2.3]{seitz-unipotent}, the $S(p)$-module
  $\res^S_{S(p)} T(d)$ is projective of dimension $2p$. If $U(p) <
  S(p)$ is the subgroup generated by some element of order $p$, it
  follows that $\res^S_{U(p)} T(d)$ is  free of
  rank 2 over $kU(p)$. Now (1) is immediate, since $u$ is $S$-conjugate to an
  element in $U(p)$.
  
  The Lie algebra case is essentially the same argument. Write $\slie
  = \Lie(S)$.  One then knows that $T(d)$ is a direct summand of
  $L(r+1) \tensor L(p-1)$, where $d = r+p$.  Since the Steinberg
  module $L(p-1)$ is projective as a module for the restricted
  enveloping algebra of $\slie$ \cite[Prop.  II.10.2]{JRAG}, the
  same holds for $T(d)$. Put $\ulie = kX \subset \slie$; together with
  \cite[Cor.  3.4]{jantzen98:Lie-prime}, the above shows that $T(d)$
  is free of rank two over $\Lambda$, the restricted enveloping
  algebra of $\ulie$. Since $\Lambda$ is a truncated polynomial
  algebra, $X$ must act with the indicated partition.
\end{proof}

Let $u,X$ as in the previous proposition. On a simple module 
$L(d) = T(d)$ with $d<p$,  both $u$ and $X$ act as a single Jordan block.
We thus obtain:

\begin{cor}
  \label{cor:tilting-jordan-blocks}
  Let $(\rho,V)$ be a tilting module for $S$ such that each weight $d$
  of a maximal torus of $S$ on $V$ satisfies $d \le 2p-2$. If $1 \not
  = u \in S$ is unipotent and $0 \not = X \in \Lie(S)$ is nilpotent,
  then the partition of $\rho(u)$ is the same as the partition of
  $d\rho(X)$.
\end{cor}

\subsection{Some results on classical groups}

Let $\Omega = \Omega(V)$ be a classical group in good characteristic
(\S\ref{sub:classical}).
\begin{prop}
  \label{prop:classical-a1-conj-class}
  Let $u \in \Omega$ be a unipotent element of order $p$, and let
  $\lambda = (\lambda_i)_{1 \le i \le t}$ be the partition of $u$.
  There is a unique $\Omega$-conjugacy class of homomorphisms
  $\psi:\SL_2(k) \to \Omega$ with the following properties:
  \begin{enumerate}
  \item The image of $\psi$ meets the conjugacy class of $u$.
  \item The character of the $\SL_2(k)$-module $(\nu \circ \psi,V)$ is
    \begin{equation*}
      \chi(\lambda_1 - 1) + \chi(\lambda_2 - 1) + \cdots + \chi(\lambda_t - 1).
    \end{equation*}
  \end{enumerate}
\end{prop}

\begin{proof}
  For existence of such a $\psi$, see \cite[Prop.
  4.1]{seitz-unipotent}. To see that any two such homomorphisms are
  conjugate, see the proof of \cite[Prop. 7.1]{seitz-unipotent}.
\end{proof}

\begin{prop}
  \label{prop:classical-tilting}
  Let $\psi:\SL_2(k) \to \Omega$ be a homomorphism satisfying (1) and
  (2) of Proposition \ref{prop:classical-a1-conj-class}.  Then $(\Ad
  \circ \psi,\olie)$ is a tilting module for $\SL_2(k)$.
  Moreover,  each weight $d$ of this representation satisfies $d \le 2p-2$.
\end{prop}

\begin{proof}
  Since the tensor product of tilting modules is again a tilting
  module \cite[Prop. 1.2]{Donk3}, and since each weight $e$ of $(\nu
  \circ \psi,V)$ satisfies $e \le p-1$, Lemmas
  \ref{lem:classical-adjoint-modules} and \ref{lem:tensor-square}
  yield the result.
\end{proof}

Let $X \in \olie$ be nilpotent, and let $\lambda_X$ be the partition
of $X$ (recall that for classical groups, the nilpotent and unipotent
classes are classified by their partitions; see e.g.
\cite[\S7.11]{hum-conjugacy} or \cite[Theorem 1.6]{jantzen:Nilpotent}
for a description of those $\lambda$ which may arise as $\lambda_X$).
As in \S \ref{sub:distinguished}, let $\phi$ be a co-character
associated with $X$.  In fact, the existence of such a $\phi$ is quite
easy to prove for a classical group; see \cite[\S
3.5]{jantzen:Nilpotent}.

\begin{prop}
  \label{prop:classical-associate-cochar}
  Let $X \in \olie$ be nilpotent with partition $\lambda =
  (\lambda_i)$. Suppose that $\lambda_i \le p$ for each $i$.  Let
  $\psi:\SL_2(k) \to \Omega$ be a homomorphism such that $X$ is in the
  image of $d\psi$, and such that the restriction of $\psi$ to a
  maximal torus of $\SL_2(k)$ is a co-character associated to $X$.
  Then $\psi$ satisfies (1) and (2) of Proposition
  \ref{prop:classical-a1-conj-class} for the unipotent class with
  partition $\lambda$.
\end{prop}

\begin{proof}
  The explicit recipe \cite[\S 3.5]{jantzen:Nilpotent} for a
  co-character associated to $X$ implies that the character of the
  $\SL_2(k)$-representation $(\Ad \circ \psi,\olie)$ is as in (2) of
  Proposition \ref{prop:classical-a1-conj-class}. For condition (1),
  it suffices to observe that $(\Ad \circ \psi,\olie)$ is a
  restricted, semisimple representation of $\SL_2(k)$ (by the linkage
  principle), so that a unipotent element in the image of $\psi$
  indeed acts with the partition $\lambda$.
\end{proof}

\subsection{The main result for unipotent elements of order $p$.}

Let $G$ be quasisimple in very good characteristic, and let $\theta:\NN \to
\UU$ be a $G$-equivariant homeomorphism as in Proposition
\ref{prop:equivariant-iso}.  Suppose
that $u = \theta(X)$, and that $u^p=1$. According
to Proposition \ref{prop:equivariant-iso}, this is equivalent to $X^{[p]}=0$.
\begin{theorem}
  \label{theorem:main-theorem-p-nilpotent}
  The partition of $\Ad u$ is the same as that of $\ad X$.
\end{theorem}

Before giving the proof, we note the following lemmas which provide a
helpful reduction in case the root system of $G$ is classical.

\begin{lem}
  \label{lem:adjoint-group}
  Let $\pi:G \to G_{\text{ad}}$ be the isogeny to the adjoint group.
  The partition of $\Ad(u)$
  is the same as that of $\Ad(\pi(u))$, and the partition of $\ad(X)$ is
  the same as that of $\ad(d\pi(X))$.
\end{lem}

\begin{proof}
  Since $p$ is very good, $d\pi$ is a bijection and the
  result is clear.
\end{proof}

\begin{lem}
  \label{lem:may-use-classical-group}
  Suppose that the root system $R$ of $G$ is classical ($A,B,C$ or
  $D$), that the characteristic is very good for $G$, and that $G$ is
  of adjoint type. There is a classical group $\Omega$ (\S
  \ref{sub:classical}), an $\Ad(\Omega)$ invariant decomposition
  $\olie = \olie_1 \oplus \olie_2$,  a unipotent $u' \in \Omega$, and a
  nilpotent $X' \in \olie$  such that
  \begin{enumerate}
  \item $u'$ and $X'$ correspond under a Springer isomorphism
    (Propoosition \ref{prop:equivariant-iso}) for $\Omega^o$,
  \item the partition of $\Ad(u)$ coincides with that of
    $\Ad(u')_{\mid \olie_1}$, 
  \item the partition of $\ad(X)$ coincides with that of
    $\ad(X')_{\mid \olie_1}$,
  \item $\olie_2$ is either zero, or $\Ad(\Omega)$ acts trivially on
    it.
  \end{enumerate}
\end{lem}

\begin{proof}
  By the previous lemma, we may replace $G$ be the corresponding
  adjoint group. Since $G$ has a ``classical type'' root system, it is
  well known that $G$ is a quotient of a classical group $\Omega$ by a
  central subgroup. It suffices to find a decomposition $\olie =
  \olie_1 \oplus \olie_2$ such that the quotient map $d\pi$ induces an
  $\Omega$-equivariant isomorphism $\olie_1 \iso \glie$, and such that
  $\olie_2$ is 0 or a trivial module. Indeed, \cite[Prop.
  26]{Mc:sub-principal} then shows that $d\pi$ will also induce
  $\Omega$-equivariant isomorphisms between the nilpotent variety of
  $\olie$ and that of $\glie$, and that $\pi$ will induce an
  isomorphism between the unipotent variety of $\Omega$ and that of
  $G$, and the lemma will follow.
  
  When $R=B,C$ or $D$, the map $\pi$ is an isogeny, and since $\olie =
  \olie_1$ and $\glie$ are simple representations, $d\pi$ induces an
  isomorphism on Lie algebras.
  
  When $R = A_r$, we may take $\Omega = \GL_{r+1}(k)$ and our
  hypothesis shows that $p$ does not divide $r+1$. Set $\olie_1 =
  \lie{sl}_{r+1}(k)$ and let $\olie_2$ be the space spanned by the
  identity matrix.  Then $\olie = \olie_1 \oplus \olie_2$ is an
  $\Omega$-stable decomposition, and $\Omega$ acts trivially on
  $\olie_2$.  Moreover, $\olie_1$ and $\glie$ are simple modules, so
  $d\pi$ induces an isomorphism $\olie_1 \iso \glie$, which completes
  the proof of the lemma.
\end{proof}

\begin{proof}[Proof of Theorem \ref{theorem:main-theorem-p-nilpotent}]
  Fix $u \in G$ unipotent of order $p$.  Let $\phi:\SL_2(k) \to G$ be
  a homomorphism whose image contains $u$ and is a good $A_1$
  subgroup, in the sense of \cite{seitz-unipotent}; this means that
  the weights of the $\SL_2(k)$-module $(\Ad \circ \phi,\glie)$ are
  all $\le 2p-2$ (the existence of such a $\phi$ is proved in
  \emph{loc. cit.}).  It follows from Theorem 1.1 of \emph{loc. cit.}
  that $(\Ad \circ \phi,\glie)$ is a tilting module for $\SL_2(k)$
  (recall that we suppose $p$ to be \emph{very good} so that the
  exceptional situation in \emph{loc. cit.} does not occur).
  
  Now let $\psi:\SL_2(k) \to G$ be a ``sub-principal'' homomorphism,
  as constructed in \cite{Mc:sub-principal}, whose image contains $u$.
   Write $u = \theta(X)$ for $X \in \NN$. Then \cite[Theorem
  12]{Mc:sub-principal} shows that the image of $d\psi$ meets the
  nilpotent orbit containing $X$.  So Theorem
  \ref{theorem:main-theorem-p-nilpotent} will follow from Corollary
  \ref{cor:tilting-jordan-blocks} once we show that $\psi$ is
  conjugate to $\phi$.  According to \cite[Prop.
  7.2]{seitz-unipotent}, all good $A_1$-subgroups meeting the class of
  $u$ are conjugate in $G$. So we are reduced to proving that the
  image of $\psi$ is a good $A_1$-subgroup.
  
  By \cite[Theorem 2]{Mc:sub-principal}, one knows that the weights of
  $(\Ad \circ \psi,\glie)$ are determined by the labeled Dynkin
  diagram of $u$.  When $G$ has an exceptional-type root system (E,F,
  or G) it follows from the proof of \cite[Prop.
  4.2]{seitz-unipotent} that those weights are $\le 2p-2$, so that the
  image of $\psi$ is indeed a good $A_1$-subgroup.
  
  When $G$ has a classical-type root system (A, B, C, or D), we may
  suppose, according to the previous lemma, that $G$ is a classical
  group.  By \cite[Theorem
  2]{Mc:sub-principal}, the restriction of $\psi$ to a maximal torus
  of $\SL_2(k)$ is a co-character associated to a conjugate of $X$.
  Thus Proposition \ref{prop:classical-associate-cochar} implies that
  $\psi$ satisfies the conditions of Proposition
  \ref{prop:classical-a1-conj-class} for $u$. Then Proposition
  \ref{prop:classical-tilting} shows that the image of $\psi$ is
  a good $A_1$ subgroup, as desired.
\end{proof}

\subsection{Further remarks on conjugacy of $A_1$-subgroups}

In this section, $G$ denotes a connected, quasisimple group.
If $G$ has a classical-type root system, we suppose
that the characteristic is very good for $G$. If
$G$ has an exceptional-type root system $R$, we suppose
that $p>N$  where $N=3,3,5,7,7$ when $R = G_2,F_4,E_6,E_7,E_8$
respectively.

\begin{theorem}
  Let $X \in \glie$ be a distinguished nilpotent element.  There is a
  unique $G$-conjugacy class of homomorphisms $\phi:\SL_2(k) \to G$
  such that the image of $d\phi$ meets the orbit of $X$. For any such
  $\phi$, the image of $\phi$ meets the unipotent class containing
  $\theta(X)$, where $\theta$ is an equivariant isomorphism as in
  Proposition \ref{prop:equivariant-iso}.
\end{theorem}

\begin{proof}
  That there is at least one such $\phi$ follows from \cite[Theorem
  12]{Mc:sub-principal}; moreover it follows from \emph{loc. cit.}
  that the image of $\phi$ meets the class of $\theta(X)$, so the last
  assertion of the theorem will follow from conjugacy.  If $\phi$ is
  any such homomorphism, we may suppose that $X = d\phi
  \begin{pmatrix}
    0 & 1 \\
    0 & 0
  \end{pmatrix}$. Since $X$ is distinguished, it 
  follows from Lemma \ref{lem:sl2-disting} that $\nu(t) =
  \phi(\operatorname{diag}(t,1/t)):k^\times \to G$ is a
  co-character associated to $X$.
  
  Suppose first that $G$ has a classical root system.  Since
  any central extension of $\SL_2(k)$ is split, we may replace $G$
  by the corresponding adjoint group, as in Lemma
  \ref{lem:adjoint-group}.  As in the proof of Lemma
  \ref{lem:may-use-classical-group}, we may therefore suppose that $G
  = \Omega = \Omega(V)$ is a classical group. In this case, the
  existence of at least one $\phi$ follows also from Proposition
  \ref{prop:classical-a1-conj-class} (this is more elementary than
  \cite{Mc:sub-principal}). The remaining assertions follow from
  Proposition \ref{prop:classical-associate-cochar}.
  
  Now supose that $G$ has an exceptional root system.  Improving on a
  result of Liebeck--Seitz, Lawther and Testerman prove under our
  assumptions on $p$ that $\phi$ is determined up to $G$-conjugacy by
  $\nu$; see \cite[\S2 and Theorem 4]{Lawther-Testerman}. This
  completes the proof.
\end{proof}

\begin{rem}
  There are of course homomorphisms $\phi:\SL_2(k) \to G$ whose image
  meets a distinguished unipotent class, but for which $d\phi$ fails
  to meet the corresponding nilpotent class. To give an example,
  recall for $0 \le m < p$ that the $\SL_2(k)$-module $L(m)$ carries
  an invariant non-degenerate bilinear form $\beta_n$ which is
  alternating when $m$ is odd, and symmetric when $m$ is even.  Now,
  suppose that $p > 2n > 0$, and consider the $\SL_2(k)$-module
  $W=L(2n) \tensor L(1)^{[1]}$ (where the exponent in the second
  factor denotes a Frobenius twist).  The tensor product $\beta =
  \beta_{2n} \tensor \beta_1$ is an invariant alternating form on $W$.
  If $\rho:\SL_2(k) \to \SP(W,\beta)$ is the defining homomorphism,
  then the image of $\rho$ meets the (distinguished) unipotent class
  with partition $(2n+2,2n)$ on $W$, but $X = d\rho
    \begin{pmatrix}
      0 & 1 \\
      0 & 0 
    \end{pmatrix}$
    has partition $(2n+1,2n+1)$. Note that $X$ is \emph{not}
    distinguished.
\end{rem}

\section{Formal groups and nilpotent endomorphisms}

\label{sec:similar}

Fix in this section an arbitrary field $k$.

\subsection{Similar nilpotent endomorphisms}

Denote by $\CC$ the category whose objects are pairs
$(V,\phi)$ where $V$ is a (finite dimensional) $k$-vector space and
$\phi$ is a nilpotent $k$-endomorphism; a morphism $f:(V,\phi)\to (W,\psi)$ in $\CC$ is a
linear map which intertwines $\phi$ and $\psi$.

We begin with the simple observation
\begin{lem}
  \label{lem:series-conj}
  Let $f \in tk\pow{t} \setminus t^2k\pow{t}$.  For each $(V,\phi)$ in
  $\CC$, there is an isomorphism $(V,\phi) \iso (V,f(\phi))$.
\end{lem}

\begin{proof}
  Since any subspace invariant by $\phi$ is also invariant by
  $f(\phi)$, we may suppose that $\phi$ acts as a single Jordan block.
  We identify $(V,\phi)$ with $(k[t]/(t^n),\mu_t)$ where $\mu_t$
  is multiplication by (the image of) $t$, so we must find a linear
  automorphism $\Lambda$ of $k[t]/(t^n)$ with $\Lambda \circ \mu_t =
  \mu_{f(t)} \circ \Lambda$.
  
  It is well known that $f$ has a ``compositional inverse''; i.e.  a
  series $g \in tk\pow{t}$ with $f(g(t)) = t$ and
  $g(f(t)) = t$. It follows that the algebra homomorphism
  $\Lambda:k[t]/(t^n) \to k[t]/(t^n)$ satisfying $t \mapsto f(t)$ is
  an automorphism. Since $\Lambda$ clearly intertwines $\mu_t$ and
  $\mu_{f(t)}$, the proof is complete.
\end{proof}

Now suppose that $Y_1,\dots,Y_m$ are indeterminants, and put
$A = k\pow{Y_1,\dots,Y_m}$. For an $m$-tuple of positive integers $\vec{r} =
(r_1,\dots,r_m)$, put
\begin{equation*}
  A_{\vec{r}} = A/(Y_1^{r_1},\dots,Y_m^{r_m}).
\end{equation*}
If $\xi_1,\dots,\xi_r \in k^\times$ and $f_1,\dots,f_r$ are formal
series without constant term, there is a uniquely determined
continuous algebra homomorphism $\Lambda:A \to A$ satisfying
$\Lambda(Y_i) = Y_i(\xi_i + f_i)$.  Evidently $\Lambda$ induces an
algebra homomorphism $\Lambda_{\vec r}:A_{\vec r} \to A_{\vec r}.$
\begin{prop}
  \label{prop:formal-calc}
  $\Lambda$ and $\Lambda_{\vec r}$ are algebra automorphisms.
\end{prop}

\begin{proof}
  Giving $A$ the $(Y_1,\dots,Y_m)$-adic filtration, the
  associated graded ring $\gr A$ identifies with $A$, and  $\gr
  \Lambda$ is given by $y_i \mapsto \xi_i y_i$, hence is an
  automorphism.  Since $A$ is complete and Hausdorff in the
  $(Y_1,\dots,Y_m)$-adic topology, \cite[II.A.4 Prop
  6]{serre:local-algebra} shows that $\Lambda$ is surjective. It
  follows that $\Lambda_{\vec{r}}$ is surjective;  since $A_{\vec{r}}$
  is finite dimensional, $\Lambda_{\vec{r}}$ is an
  automorphism. If $f \in \ker \Lambda$,  the image
  of $f$ in $A_{\vec{r}}$ is 0 for each $\vec{r}$; thus $f = 0$ and
  the proposition follows.
\end{proof}

\begin{theorem}
  \label{theorem:tensor-result}
  Let $F(u,v) \in k\pow{u,v}$ satisfy $F(u,v) \congruent \xi_1 u +
  \xi_2 v \pmod{(u,v)^2}$ with $\xi_1,\xi_2 \in
  k^\times$, and suppose that $(V,\phi)$ and $(W,\psi)$ are objects in
  $\CC$. Then there is an isomorphism in $\CC$
  \begin{equation*}
    (V \tensor W, \phi \tensor 1 + 1  \tensor \psi) \iso
    (V \tensor W ,F(\phi \tensor 1,1 \tensor \psi) ).
  \end{equation*}
\end{theorem}

\begin{proof}
  Write $X=\phi \tensor 1 + 1 \tensor \psi$ and $X'=F(\phi \tensor 1,1
  \tensor \psi)$; we must show that $X$ and $X'$ are conjugate by the
  adjoint action of $\GL(V \tensor W)$.  If $V = V_1 \oplus V_2$ with
  each $V_i$ invariant by $\phi$, then $V_i \tensor W$ is invariant by
  both $X$ and $X'$ for $i=1,2$.  So we may suppose that $\phi$ (and
  also $\psi$) acts as a single Jordan block.  We identify $(V,\phi)$
  with $(k[Y]/(Y^n),\mu_y)$ and $(W,\psi)$ with $(k[Z]/(Z^m),\mu_z)$,
  where we write $y$ and $z$ for the images of $Y$ and $Z$ in the
  quotient algebras.
  
  We must find a linear automorphism $\Lambda$ of $$B=k[Y]/(Y^n) \tensor
  k[Z]/(Z^m) = k[Y,Z]/(Y^n,Z^m)$$ satisfying $\Lambda \circ (\mu_{y +
  z}) = \mu_{F(y,z)} \circ \Lambda$.

  Choose formal series 
    $H_1(u,v),H_2(u,v) \in uk\pow{u,v} +  vk\pow{u,v}$
  such that $F(u,v) = \xi_1 u + \xi_2 v + H_1(u,v)u +H_2(u,v)v$ [in general,
  $H_1$ and $H_2$ are not uniquely determined by this condition].
  Then according to Proposition \ref{prop:formal-calc}, there is an
  algebra automorphism $\Lambda:B \to B$ satisfying $\Lambda(y) =
  f_1(y,z)$ and $\Lambda(z)= f_2(y,z)$ where
  \begin{equation*}
    f_1(Y,Z) = Y(\xi_1 + H_1(Y,Z)) \quad \text{and} \quad
    f_2(Y,Z) = Z(\xi_2 + H_2(Y,Z)).
  \end{equation*}
  
  Since
  \begin{equation*}
    \Lambda(y+z) = f_1(y,z) + f_2(y,z) = \xi_1 y + \xi_2 z + yH_1(y,z) + zH_2(y,z) 
    = F(y,z),
  \end{equation*}
  $\Lambda$ intertwines $\mu_{y+z}$ and $\mu_{F(y,z)}$ as desired.
\end{proof}

\subsection{Formal groups}
\label{sub:formal-groups}

A \emph{formal group law} over $k$ is a series $\formal(u,v) \in
k\pow{u,v}$ satisfying $\formal(u,v) \congruent u+v \pmod{(u,v)^2}$
and certain other axioms that are written down in a number of places,
se e.g.
 \cite{Hazewinkel-formalgroups}. Note that $\formal(u,v)$ is
automatically commutative; see \cite[Theorem
I,1.6.7]{Hazewinkel-formalgroups}.

The representation ring $R_\formal$ of $\formal$ introduced in
\cite{Fossum-formalgroups} is first of all the Abelian group
generated by symbols $[M]$ for each object $M$ in $\CC$ subject to
relations as follows: $[M]$ and $[N]$ are identified if $M \iso N$.
The addition is given by the direct sum in $\CC$:
\begin{equation*}
  [(V,\phi)] + [(W,\psi)] =   [(V \oplus W,\phi \oplus \psi)].
\end{equation*}
The tensor product
(with respect to $\formal$) of two objects in $\CC$ is defined by
\begin{equation}
  \label{eq:formal-tensor}
  (V,\phi) \tensor_\formal (W,\psi) = (V \tensor W,\formal(\phi \tensor 1,1 \tensor \psi)).
\end{equation}
The product in $R_\formal$ is defined by
$[(V,\phi)] \cdot [(W,\psi)] = [(V,\phi) \tensor_\formal (W,\psi)]$;
the  axioms of a formal group imply that $R_\formal$ is indeed a ring.

Let $J_n$ denote the class in $R_\formal$ obtained from the object
$(k^n,\phi_n)$ where $\phi_n$ is a single $n\times n$ Jordan block.
The $J_n$ form a $\Z$-basis for $R_\formal$, and the multiplication in
$R_\formal$ is completely determined by the structure constants
$a^i_{n,m}$ for the products $[J_n]\cdot[J_m] = \sum_i a^i_{n,m}
[J_i]$ for each $n,m \ge 0$.

Now suppose that $k$ has characteristic $p>0$.  For each $n \ge 1$,
the full subcategory of $\CC$ consisting of those objects $(V,\phi)$
such that $\phi^{p^n} = 0$ ($p^n$-nilpotent objects) is 
closed under the tensor product \eqref{eq:formal-tensor} for
any formal group $\formal$.

Let $C$ be a cyclic group of order $p^n$ with generator $g$. The
category of $kC$ modules is then equivalent to the $p^n$-nilpotent
subcategory; the $C$-module $V$ corresponds to the object $(V,g-1)$ of
$\CC$. The usual tensor product for group representations corresponds
to the tensor product \eqref{eq:formal-tensor} with respect to the
\emph{multiplicative} formal group law $\formal_m = u + v + uv$.

Let $\witt_n$ be the Abelian $p$-Lie algebra $\sum_{i=0}^{n-1} k
X^{[p^i]}$ (with the indicated $p$-power map). The category of
restricted $\witt_n$-representations is equivalent to the
$p^n$-nilpotent subcategory. The usual tensor product of Lie algebra
representations corresponds to the tensor product
\eqref{eq:formal-tensor} with respect to the additive formal group law
$\formal_a = u + v$.

\begin{cor}[Fossum]
  \label{cor:rep-ring}
  Let $\formal_1$ and $\formal_2$ be any two formal group laws over
  $k$.  The structure constants for the representation rings
  $R_{\formal_1}$ and $R_{\formal_2}$ are the same.
\end{cor}

\begin{proof}
  This is an immediate consequence of Theorem
  \ref{theorem:tensor-result}.
\end{proof}

In characteristic 0, all 1-dimensional formal group laws are
isomorphic (see \cite[Theorem 1.6.2]{Hazewinkel-formalgroups}) so one
has in that case a simple proof of the corollary.  The structure
constants in characteristic 0 are the ``Clebsch-Gordan'' coefficients;
see \cite[(2.3)]{Srinivasan-Modular}.  In characteristic $p>0$, there
are non-isomorphic formal group laws; for instance, $\formal_a \not \iso \formal_m$.

The representation ring of the multiplicative law $\formal_m$ has been
studied by J. A.  Green in \cite{Green-RepRing} and B. Srinivasan in
\cite{Srinivasan-Modular}; see also the references in
\cite{Fossum-formalgroups}. Green computed enough of the structure
constants to show that $R_{\formal_m} \tensor_\Z \C$ is a semisimple
algebra. Srinivasan determined the structure constants
of $R_{\formal_m}$ explicitly.

The above corollary was proved by Fossum in \cite[III]{Fossum-formalgroups}.
He showed that the eight formulas obtained by Green in
\cite{Green-RepRing} hold for any formal group law; the corollary 
follows since those formulas determine the multiplication in
$R_\formal$ uniquely.  This proof is in the spirit of ``modular
representation theory''; it exploits the fact that $(k^q,\phi_q)$ is a
projective object in a suitable subcategory of $\CC$ when $q = p^m$.

\subsection{A result on symmetric series}

Let $\Sigma_m$ be the symmetric group on $m$ letters.  Fix
$Y_1,\dots,Y_m$ indeterminants (each given degree 1), and consider the
graded algebras $A = k[Y_1,\dots,Y_m]$ and $\hat A =
k\pow{Y_1,\dots,Y_m}$.

The rule $\sigma(Y_i)= Y_{\sigma^{-1}(i)}$ for $\sigma \in \Sigma_m$
defines representations of $\Sigma_m$ by graded algebra automorphisms
on $A$ and by continuous graded algebra automorphisms on $\hat A$. We
use the notations $A^{\Sigma_m}$ and $\hat{A}^{\Sigma_m}$ for the
subalgebras of invariants. $A_d$ and $\hat{A}_d$ denote the respective
$d$-th homogeneous components, $A_{\ge d}$ is $\bigoplus_{e\ge d} A_e$,
and $\hat{A}_{\ge d}$ is the ideal of $\hat{A}$ generated by $A_{\ge
  d}$ (the $d$-th power of the unique maximal ideal of $\hat{A}$).

\begin{lem}
  Suppose that $m!$ is invertible in $k$.  Let $H \in A_d^{\Sigma_m}$
  be a homogeneous invariant. Then there are elements $H_i \in A_d$
  for $1 \le i \le m$ such that $H = \sum_i H_i$, $Y_i$ divides $H_i$,
  and $\sigma H_i = H_{\sigma^{-1}(i)}$ for each $\sigma \in
  \Sigma_m$.
\end{lem}

\begin{proof}
  Observe that if $H',H'' \in A_d^{\Sigma_m}$ and the conclusion of
  the lemma holds for $H'$ and $H''$, then it holds for $H' + H''$.
  Similarly, if $H' \in A_d^{\Sigma_m}$ and the conclusion of the lemma
  holds for $H'$, then it holds for $H' \cdot H''$ for any $H'' \in
  A_e^{\Sigma_m}$.
  
  Since $A^{\Sigma_m}$ is a polynomial algebra in the elementary
  symmetric polynomials $s_j \in A_j$, $1 \le j \le m$, we may
  therefore suppose $H=s_j$. Thus we fix $1 \le j \le m$.

  For  $I \subseteq \{1,\cdots,m\}$, let $Y_I = \prod_{k \in I} Y_k$.
  For each $1 \le i \le m$, we define
  \begin{equation*}
    H_i = (1/j)\sum_{|I|=j,\ i \in I} Y_I;
  \end{equation*}
  note that $j$ is invertible in $k$ by assumption.  Evidently $Y_i$
  divides $H_i$.  It is a simple matter to verify that $\sum_{i=1}^m
  H_i = s_j.$ Let $\sigma \in \Sigma_m$. Since $\sigma Y_I = \prod_{k
    \in I} Y_{\sigma^{-1}(k)}$, it follows that $\sigma H_i =
  H_{\sigma^{-1}(i)}$.
\end{proof}

\begin{prop}
  \label{prop:symmetric-series}
  Suppose that $m!$ is invertible in $k$ and that $f \in
  \hat{A}^{\Sigma_m}$ satisfies $f \congruent Y_1 + \cdots + Y_m
  \pmod{\hat{A}_{\ge 2}}$ with $\xi \in k^\times$. Then we may find series
  $f_i \in \hat{A}$ for $1 \le i \le m$ such that
  \begin{enumerate}
  \item   $f_i \congruent Y_i \pmod{\hat{A}_{\ge 2}}$,
  \item $Y_i \mid f_i$,
  \item $\sigma f_i = f_{\sigma^{-1}(i)}$ for all $\sigma \in \Sigma_m$.
  \item $f = f_1 + \cdots + f_m$.
  \end{enumerate}
\end{prop}

\begin{proof}
  Since $\Sigma_m$ acts by graded algebra automorphisms, $f$ is
  invariant if and only if its homogeneous components are. So the
  proposition follows from the previous lemma.
\end{proof}

\subsection{Exterior and symmetric powers in $\CC_\formal$.}
\label{sub:multilinear-result}
Let $\formal(u,v)$ be a formal group law as before.  Let
$Y_1,Y_2,\dots$ be indeterminants, and put
\begin{equation*}
  \tensorpoly^1_\formal(Y_1) = Y_1, \quad
  \tensorpoly^m_\formal(Y_1,\dots,Y_m) = 
  \formal(\tensorpoly^{m-1}_\formal(Y_1,\dots,Y_{m-1}),Y_m), \quad m \ge 2.
\end{equation*}
Thus $\formal(Y_1,Y_2) = \tensorpoly^2_\formal(Y_1,Y_2)$,, and if
$(V,\phi)$ is an object of $\CC$, the $m$-fold power of $(V,\phi)$ for
the tensor product \eqref{eq:formal-tensor} is $(V^{\tensor
  m},\tensorpoly^m_\formal(\phi))$, where
$\tensorpoly^m_\formal(\phi)$ is obtained from
$\tensorpoly^m_\formal(Y_1,\dots,Y_m)$ by specializing $Y_i \mapsto
1^{\tensor (i-1)} \tensor \phi \tensor 1^{\tensor (m-i)}$.

Recall that $\formal_a(u,v) = u+v$ denotes the additive formal group
law. We have for each $m \ge 1$:
\begin{equation*}
  \tensorpoly^m_{\formal_a}(Y_1,\dots,Y_m) = Y_1 + \cdots + Y_m.
\end{equation*}
For any $\formal(u,v)$ we have
\begin{equation*}
  \tensorpoly^m_{\formal}(Y_1,\dots,Y_m) \congruent Y_1 + \cdots + Y_m 
  \pmod{k\pow{Y_1,\dots,Y_m}_{\ge 2}}.
\end{equation*}

There is a linear representation of the symmetric group $\Sigma_m$ on
$V^{\tensor m}$ given by
\begin{equation}
  \label{eq:symmetric-group-tensor-action}
  \sigma(v_1\tensor\cdots \tensor v_m) =
  v_{\sigma^{-1}(1)} \tensor \cdots \tensor v_{\sigma^{-1}(m)}
\end{equation}
for $\sigma \in \Sigma_m$ and $v_i \in V$, $1 \le i \le m$. Any map $X
\in \End_{\Sigma_m}(V^{\tensor m})$ induces maps $X' \in
\End(\bwedge^m V)$ and $X'' \in
\End(\Sym^m V)$. 

Since any 1-dimensional formal group law $\formal(u,v)$ is
commutative, we have
\begin{equation*}
  \tensorpoly^m_\formal(Y_1,\dots,Y_m) \in k\pow{Y_1,\dots,Y_m}^{\Sigma_m}.
\end{equation*}
It follows that $\tensorpoly^m_\formal(\phi)$ is a
$\Sigma_m$-homomorphism and hence induces an endomorphism
$\bwedge^m_\formal(\phi)$ of the exterior power $\bwedge^m V$ and an
endomorphism $\Sym^m_\formal(\phi)$ of the symmetric power $\Sym^m V$.

\begin{theorem}
  \label{theorem:multilinear-result}
  Let $\formal$ be a formal group law, and   $(V,\phi)$ an object in
  $\CC$.  Assume that $m \ge 1$ is an integer such that $m!$
  is non-zero in $k$.  Then there are isomorphisms in $\CC$
  \begin{equation*}
    (\bwedge^m V,\bwedge^m_{\formal_a}(\phi))
    \iso
    (\bwedge^m V,\bwedge^m_{\formal}(\phi))
  \end{equation*}
  and
  \begin{equation*}
    (\Sym^m V, \Sym^m_{\formal_a}(\phi)) 
    \iso
    (\Sym^m V, \Sym^m_{\formal}(\phi)),
  \end{equation*}
  where $\formal_a = u + v$ is the additive formal group law.
\end{theorem}

\begin{proof}
  If $V = W \oplus W'$ where $W$ and $W'$ are proper subspaces
  invariant by $\phi$, one has the decomposition
      \begin{equation*}
    \bwedge^m V = 
    \displaystyle \bigoplus_{i+j=m}  
      \bwedge^i W \tensor  \bwedge^j W'
  \end{equation*}
  into subspaces invariant by both $\bwedge^m_\formal(\phi)$ and
  $\bwedge^m_{\formal_a}(\phi)$.

  Let $W_{i,j} = \bwedge^i W \tensor \bwedge^j W'$ with $i+j = m$. On
  $W_{i,j}$, $\bwedge^m_\formal(\phi)$ acts as
  $$\formal(\bwedge^i_\formal(\phi) \tensor 1, 1 \tensor
  \bwedge^j_\formal(\phi)),$$
  which by Theorem
  \ref{theorem:tensor-result} is similar to $\bwedge^i_\formal(\phi)
  \tensor 1 +1 \tensor \bwedge^j_\formal(\phi)$. If the result were
  known for $\bwedge^i W$ and $\bwedge^j W'$, one would know this last
  endomorphism to be similar to $\bwedge^i_{\formal_a}(\phi) \tensor 1
  +1 \tensor \bwedge^j_{\formal_a}(\phi)$, which is precisely the
  restriction of $\bwedge^m_{\formal_a}(\phi)$ to $W_{i,j}$.  We
  deduce that if the theorem is known for $W$ and $W'$, then it holds
  for $V$. Similar remarks hold for the symmetric powers.  Thus, we
  may suppose that $\phi$ acts as a single Jordan block on $V$.
  
  Let $B = k[Y_1,\dots,Y_m]/(Y_1^n,\dots,Y_m^n)$ where $n = \dim V$;
  writing $y_i$ for the image of $Y_i$ in $B$, we identify
  $(V^{\tensor m},\tensorpoly^m_\formal(\phi))$ with $(B,\mu)$ where
  $\mu$ is multiplication by $\tensorpoly^m_\formal(y_1,\dots,y_m)$.
  Similarly, we identify $(V^{\tensor m},\tensorpoly^m_\formal(\phi))$
  with $(B,\mu')$ where $\mu'$ is multiplication by $y_1 + \cdots +
  y_m$.  These identification are isomorphisms of $\Sigma_m$
  representations, where $\Sigma_m$ acts on $V^{\tensor m}$ as in
  \eqref{eq:symmetric-group-tensor-action}, and on $B$ by the action
  induced from that on the polynomial algebra.
  
  The theorem will follow if we find a $k\Sigma_m$-linear automorphism
  $\Lambda$ of $B$ such that $\Lambda \circ \mu = \mu' \circ \Lambda$.
  Since $\tensorpoly^m_\formal(Y_1,\dots,Y_m) \in
  k\pow{Y_1,\dots,Y_m}$ satisfies the hypothesis of Proposition
  \ref{prop:symmetric-series}, we may find power series
  $f_1,\dots,f_m$ satisfying conditions 1-4 of that Proposition.  In
  view of conditions 1 and 2, Proposition \ref{prop:formal-calc}
  implies that the rule $\Lambda(y_i) = f_i(y_1,\dots,y_m)$ defines an
  algebra automorphism of $B$.  Property 3 implies that $\Lambda$
  intertwines the $\Sigma_m$ action.  Finally, property 4 implies that
  $\Lambda$ intertwines $\mu$ and $\mu'$.  This completes the proof of
  the theorem.
\end{proof}

\begin{example}
  \label{rem:fossum-example}
  In \cite{Fossum-formalgroups}, it is proved that the classes of
  \begin{equation*}
      (\bwedge^m V, \bwedge^m_\formal(\phi)) \quad \text{and} \quad (\Sym^m V,
      \Sym^m_\formal(\phi))
  \end{equation*}
  in the representation ring are independent of $\formal$ provided
  $\phi^p = 0$, if $p$ is the characteristic of $k$.  Theorem
  \ref{theorem:multilinear-result} shows that these classes are
  independent of $\formal$ provided that $p>m$ (or $p=0$ of course).
  
  To see that some hypothesis on the characteristic is essential,
  suppose that $p=2$. Recall that we write $(V_n,\phi_n)$ for the $n$
  dimensional Jordan block in $\CC$.  Then some computer calculations
  yield:
  \begin{equation*}
    \begin{array}[t]{ll|ccc}
      n & \formal &  [V_n^{\tensor 2}, \tensorpoly_\formal^2(\phi_n)] 
        & [(\bwedge^2 V_n,\bwedge^2_\formal(\phi_n))] 
        & [(\Sym^2 V_n,\Sym^2_\formal(\phi_n))]\\
      \hline \hline
      4 & \formal_m & 4J_4 & J_2 + J_4 & 2J_4 + J_2  \\
      4 & \formal_a & 4J_4 & 2J_3 & 2J_4 + 2J_1  \\
      5 & \formal_m & 2J_8 + 2J_4 + J_1 & J_7 + J_3 & J_8 + J_4 + J_3 \\
      5 & \formal_a & 2J_8 + 2J_4 + J_1 & J_7 + J_3 & J_8 + J_4 + 3J_1 \\
      6 & \formal_m & 4J_8 + 2J_2 & J_8 + J_6 + J_1 & 2J_8 + J_4 + J_1 \\
      6 & \formal_a & 4J_8 + 2J_2 & 2J_7 + J_1 & 2J_8 + J_2 + 3J_1\\
      7 & \formal_m & 6J_8 + J_1 & 2J_8 + J_5 & 3J_8 + J_4 \\
      7 & \formal_a & 6J_8 + J_1 & 3J_7 & 3J_8 + 4J_1 \\
    \end{array}
  \end{equation*}

\end{example}

\section{Adjoint Jordan blocks for classical groups}

Let $\Omega = \Omega(V)$ be a classical group 
in good characteristic (see \S \ref{sub:classical}).

\begin{prop}
  There is a series $\e(t) \in tk\pow{t} \setminus
  t^2k\pow{t}$ such that $X \mapsto 1 + \e(X)$ defines
  a $\Omega$-equivariant isomorphism of varieties $\NN \to \UU$.
\end{prop}

\begin{proof}
  See \cite[\S6.20]{hum-conjugacy}. 
\end{proof}

Note that we identify $\NN$ and $\UU$ with subvarieties of $\End_k(V)$.
In case CG1, \emph{any} series $\e(t) \in tk\pow{t} \setminus
t^2k\pow{t}$ defines an isomorphism as in the proposition. In cases
CG2, CG3, one may take the Cayley transform series
\begin{equation*}
  \e(t) = (1-t)(1+t)^{-1} -1 = 
  2 \sum_{i=1}^\infty (-1)^{i} t^i.
\end{equation*}
We may also use in these cases the Artin-Hasse exponential series
\cite[V\S16]{SerreAG&CF}.  This follows from
\cite[7.4]{mcninch-math.RT/0007056}; see also \cite{Proud-Witt}.
 
\begin{theorem}  
  \label{theorem:classical-adjoint-partitions}
  Let $\theta:\NN \to \UU$ be an $\Omega$-equivariant
  isomorphism of varieties. Then $\Ad(\theta(X))$ and $\ad(X)$ have
  the same partition.
\end{theorem}

\begin{proof}
  We may suppose that $\theta$ is determined by a formal series $\e$
  as in the previous proposition.  Using the identifications of the
  adjoint modules in Lemma \ref{lem:classical-adjoint-modules}, the action of
  $\ad(X)$ is determined by the additive formal group law $\formal_a$.
  We have:
  \begin{equation*}
    \ad(X) = \left \{ 
      \begin{matrix}
        X \tensor 1 + 1 \tensor X^\vee & \text{in case CG1} \\[3pt]
        \bwedge^2_{\formal_a} (X) & \text{in case CG2} \\[3pt]
        \Sym^2_{\formal_a} (X) & \text{in case CG3.} \\
      \end{matrix}
    \right .
  \end{equation*}

   Similarly, we have
     \begin{equation*}
       \Ad(\theta(X)) = \left \{ 
      \begin{matrix}
        (1 + \e(X)) \tensor (1 + \e(X^\vee)) & \text{in case CG1} \\[3pt]
        1 + \bwedge^2_{\formal_m} (\e(X)) & \text{in case CG2} \\[3pt]
        1 + \Sym^2_{\formal_m} (\e(X)) & \text{in case CG3.} \\
      \end{matrix}
    \right .
  \end{equation*}
  By Lemma \ref{lem:series-conj}, $\Ad(\theta(X))$ is similar
  to 
  \begin{equation*}
      \begin{matrix}
        (1 + X) \tensor (1 + X^\vee) & \text{in case CG1,} \\[3pt]
        1 + \bwedge^2_{\formal_m} (X) & \text{in case CG2, and} \\[3pt]
        1 + \Sym^2_{\formal_m} (X) & \text{in case CG3.} \\
      \end{matrix}
  \end{equation*}
   The result now follows from Theorem \ref{theorem:tensor-result} and
   Theorem \ref{theorem:multilinear-result}.
\end{proof}

\begin{rem}
  In bad characteristic, the adjoint partition of the regular
  nilpotent class can already differ from that of the regular
  unipotent class.  Let $p=2$ and let $\Omega$ be the symplectic group
  $\SP_{2n}(k)$. If $u \in G$ is a regular unipotent element, and $X
  \in \olie$ is a regular nilpotent element, them from example
  \ref{rem:fossum-example} one obtains the following partitions:
  \begin{equation*}
    \begin{array}[t]{lll}
      n & \text{partition of $\Ad(u)$} & \text{partition of $\ad(X)$} \\
      \hline \hline
      2 & (4^2,2) & (4^2,1^2) \\
      3 & (8^3,4) & (8^2,2,1^3) \\
    \end{array}
  \end{equation*}
  Note that the dimensions of $\ker (\Ad(u) -1)$ and $\ker \ad(X)$
  agree with the tables in \cite[p. 179]{MR82d:14030}.
  
  If again $p=2$, $\Omega$ is the orthogonal group $\operatorname{O}_7(k)$,
  and $u,X$ are regular, then the same example shows that the
  partition of $\Ad(u)$ is $(8^2,5)$ while that of $\ad(X)$ is
  $(7^3)$.
\end{rem}

\section{Adjoint Jordan blocks for $G_2$}
\label{sub:g2}

In this section, we consider a group $G$ with root system of type
$G_2$. Denote by $\varpi_1,\varpi_2$ the fundamental dominant weights
in the character group of a fixed maximal torus $T$ of $G$ with respect to
a choice of simple roots $\alpha_1,\alpha_2$, where $\alpha_1$ is
\emph{short}. Suppose that $p>3$.

Since $p \not = 2$, the Weyl module $V=V(\varpi_1)$ for $G$ is
irreducible of dimension 7.  Moreover, $G$ leaves invariant a
non-degenerate quadratic form on $V$; this defines a closed embedding
$\rho:G \to H=SO(V)$. Recall that $H$ has a root system of type $B_3$.

\subsection{Simple root vectors and root subgroups}
\label{sub:simple-roots-g2}

Fix non-zero elements $x_{\alpha_i} \in \glie_{\alpha_i}$, $i=1,2$.
\begin{lem}
  \label{lem:g2-nilp-lemma}
  For suitable choices of a maximal torus of $H$ containing $T$ and
  a system $\gamma_i$ ($i=1,2,3$)  of simple roots for $H$, we have
  \begin{enumerate}
  \item $d\rho(x_{\alpha_1}) = y_{\gamma_1} + y_{\gamma_3} \in
\hlie_{\gamma_1} + \hlie_{\gamma_3}$
\item $d\rho(x_{\alpha_2}) = y_{\gamma_2} \in \hlie_{\gamma_2}$
  \end{enumerate}
  with $y_{\gamma_i} \not = 0$ for $i=1,2,3$.
\end{lem}

\begin{proof}
  The weights of the $G$-module $V$ are $\pm f_1, \pm f_2, \pm f_3$
  and $0$, where $f_1= \varpi_1$, $f_2 = \varpi_1 - \alpha_1$, and
  $f_3 = \varpi_1 - \alpha_1 - \alpha_2$.
  
  Since $V$ is a simple, restricted $G$-module, a theorem of Curtis
  implies that it is also simple as a $\glie$-module.  Let $\nlie$ be
  the Lie algebra of the unipotent radical of the Borel sub-group of
  $G$ corresponding to our choice of positive roots.  Since $V$ is
  simple for $\glie$, we must have $V^{\nlie} = V_{f_1}$.  Since
  $\nlie$ is generated as a Lie algebra by $x_{\alpha_1}$ and
  $x_{\alpha_2}$, it follows that for any weight $\gamma \not = f_1$
  of $V$, the restriction of either $x_{\alpha_1}$ or $x_{\alpha_2}$
  to $V_\lambda$ is non-zero.

  Thus, we see that 
  \begin{equation}
    \label{eq:x-a-1}
    x_{\alpha_1} V_{f_2} = V_{f_1}, \quad x_{\alpha_1} V_0 = V_{f_3},
    \quad x_{\alpha_1} V_{-f_3} = V_{0},     
    \quad x_{\alpha_1} V_{-f_1} = V_{-f_2}
  \end{equation}
  and $x_{\alpha_1}$ acts trivially on all other weight spaces of $V$.
  Similarly,
  \begin{equation}
        \label{eq:x-a-2}
    x_{\alpha_2} V_{f_3} = V_{f_2}, \quad x_{\alpha_2} V_{-f_2} = V_{-f_3}.
  \end{equation}
  and $x_{\alpha_2}$ acts trivially on all other weight spaces of $V$.
  
  Choose a non-zero weight vector $e_{\pm i} \in V_{\pm f_i}$ and $e_0
  \in V_0$; this yields a basis of $V$. Since $T$-weight spaces
  $V_\lambda$ and $V_\mu$ are orthogonal under the invariant form
  unless $\lambda + \mu = 0$, it follows from \cite[\S 23.4]{Bor1}
  that the subgroup $S$ of $H$ which acts diagonally in this basis is
  a maximal torus (containing $T$). If $\e_i \in X^*(S)$ is the
  $S$-weight of $e_i$, then the $\e_1,\e_2,\e_3$ form a $\Z$-basis for
  $X^*(S)$. Moreover, $\gamma_1=\e_1 - \e_2, \gamma_2=\e_2 - \e_3,
  \gamma_3=\e_3$ is a system of simple roots for $H$. The roots of $H$
  are: $R_H = \{\pm(\e_i \pm \e_j) \mid 1 \le i < j \le 3\} \cup\{\pm
  e_3\}$.
  
  Write $x_{\alpha_i} = \sum_{\alpha \in R_H} z_{i,\alpha} \in
  \bigoplus_{\alpha \in R_H} \hlie_\alpha$ for $i=1,2$ (since
  $x_{\alpha_i}$ is nilpotent, it has no component in $\hlie_0$).  To
  prove (1), it suffices to show that $z_{1,\alpha} = 0$ whenver
  $\alpha \not = \gamma_1,\gamma_3$; this follows from
  \eqref{eq:x-a-1}.  To prove (2), it suffices to show that
  $z_{2,\alpha} = 0$ whenver $\alpha \not = \gamma_2$; this follows
  from \eqref{eq:x-a-2}.
\end{proof}

Let $U_{\alpha_i} \le G$ ($i=1,2$) be the simple root subgroups.  We
claim that $U_{\alpha_i}$ is the image of the map $\XX_{\alpha_i}:\G_a
\to G$ satisfying $\rho(\XX_{\alpha_i}(s)) = \exp(s d\rho(x_{\alpha_i}))$.  One
can see the claim in several ways. One is to realize $G$ as a
Chevalley group (see \cite{Steinberg}). Alternatively, one may realize
$G$ as the automorphisms of an 8 dimensional Cayley algebra
$\mathbf{O}$; the module $V$ is the space orthogonal to the identity
element of $\mathbf{O}$. Since the $x_{\alpha_i}$ act as derivations
of $\mathbf{O}$, their exponentials (which are defined thanks to our
assumption on $p$) yield automorphisms of $\mathbf{O}$. The images of
these exponentials yield 1 dimensional unipotent subgroups of $G$
normalized by $T$ which are then necessarily the desired root subgroups.

Since $p \not = 2$, it is well-known that the simple root subgroups
$V_{\gamma_i} \le H$ are the images of the maps $\YY_{\gamma_i}:\G_a
\to H$ given by $\YY_{\gamma_i}(s) = \exp(s y_{\gamma_i})$.
  
\begin{lem}
    \label{lem:g2-unip-lemma}
  We have $\rho(\XX_{\alpha_1}(s)) = \YY_{\gamma_1}(s)\YY_{\gamma_3}(s)$,
  and $\rho(\XX_{\alpha_2}(s)) = \YY_{\gamma_2}(s)$.
\end{lem}

\begin{proof}
  Since $[y_{\gamma_1},y_{\gamma_3}]=0$, we have $\exp(s
  d\rho(x_{\alpha_1})) = \exp(sy_{\gamma_1})\exp(sy_{\gamma_3})$ and
  the formula for $\alpha_1$ follows. The formula for $\alpha_2$ is
  simpler.
\end{proof}

\begin{prop}
  \label{prop:g2-regular-elements}
  Let $X \in \glie$ be regular nilpotent, and $u \in G$ be regular
  unipotent. Then $X$ and $u$ both act as a single Jordan block on
  $V=V(\varpi_1)$.
\end{prop}

\begin{proof}
  Recall that the regular nilpotent and regular unipotent classes for
  $H$ act as a single Jordan block on the natural module of $H$.  Also
  recall that for any reductive group $L$, if $z_\alpha \in
  \glie_\alpha$ are non-zero root vectors for a system of simple
  roots, then $\sum_\alpha z_\alpha$ is regular nilpotent.  If the
  images of $\ZZ_\alpha:\G_a \to L$ are the simple root subgroups and
  $a_\alpha \in k^\times$, then $\prod_\alpha \ZZ_\alpha(a_\alpha)$ is
  regular unipotent (one may take the product in any fixed order).
  
  In particular, $x_\alpha + x_\beta$ is regular nilpotent in $\glie$,
  and by Lemma \ref{lem:g2-nilp-lemma} we see that $d\rho(x_\alpha +
  x_\beta) = y_{\gamma_1} + y_{\gamma_2} + y_{\gamma_3}$ is regular
  nilpotent in $\hlie$. Thus $X$ acts as a single Jordan block as claimed.
  
  Similarly, $v=\XX_{\alpha_1}(1)\XX_{\alpha_2}(1)$ is regular
  unipotent in $G$, and by Lemma \ref{lem:g2-unip-lemma}, 
  \begin{equation*}
      \rho(v) = \YY_{\gamma_1}(1)\YY_{\gamma_3}(1)\YY_{\gamma_2}(1)
  \end{equation*}
  is regular unipotent in $H$. Thus $u$ also acts as a single Jordan
  block.
\end{proof}

\subsection{The adjoint representation}

Since $p \not = 3$, the Weyl module $V(\varpi_2)$ is irreducible; it is
isomorphic with the adjoint representation of $G$.

If $W$ is a $G$-module, write $\ch W = \sum_\mu \dim W_\mu \cdot
e^\mu$ for its character. For a dominant weight $\mu$, let $\chi(\mu)
= \ch V(\mu)$, where $V(\mu)$ is the Weyl module with high weight
$\mu$. Brauer's formula \cite[ex. 24.4(9)]{Hu1} yields $$\ch
V(\varpi_1) \tensor V(\varpi_1) = \chi(0) + \chi(\varpi_1)
+\chi(\varpi_2) + \chi(2\varpi_1).$$
Since $\dim V(2\varpi_1) = 27$
(by Weyl's formula) the character of $\Sym^2 V(\varpi_1)$ must be
$\chi(0) + \chi(2\varpi_1)$, and Lemma \ref{lem:tensor-square} yields
$\ch \bwedge^2 V(\varpi_1) = \chi(\varpi_1) + \chi(\varpi_2)$ (despite
the application of Lemma \ref{lem:tensor-square}, these character
formulas are valid in all characteristics).  Since $p>3$, the Weyl
modules $V(\varpi_1)$, $V(\varpi_2)$ are simple and we deduce from
\cite[Prop. II.2.14]{JRAG} that $0 = \ext^1_G(V(\varpi_1),V(\varpi_2))
= \ext^1_G(V(\varpi_2),V(\varpi_1)).$ Thus:

\begin{prop}
  \label{prop:g2-direct-sum}
   $\bwedge^2 V(\varpi_1) \iso V(\varpi_1) \oplus
  V(\varpi_2)$ is semisimple.
\end{prop}

\begin{theorem}
  \label{theorem:g2-adjoint-partitions}
  Let $G$ be a group with root system $G_2$ and suppose $p>3$.  Let
  $\theta:\NN \to \UU$ be a $G$-equivariant isomorphism.  For each $X
  \in \NN$, $\ad X$ and $\Ad \theta X$ have the same partition.
\end{theorem}

\begin{proof}
  Since $p>3$, the only nilpotent class which is not
  $p$-nilpotent is the regular class when $p=5$; see for example
  \cite[Example \S6.2]{mcninch-math.RT/0007056}.  Let $X \in \NN$ be
  regular nilpotent; hence also $u=\theta(X)$ is regular unipotent.  Let
  $\rho:G \to H$ be as before. Proposition
  \ref{prop:g2-regular-elements} implies that $d\rho(X)$ and
  $\rho(u)$ are respectively regular nilpotent and regular
  unipotent in $H$.  Thus Theorem \ref{theorem:classical-adjoint-partitions}
  and Lemma \ref{lem:classical-adjoint-modules} together imply that the
  partition of $X$ acting on $\bwedge^2 V(\varpi_1)$ is the same as
  that of $u$.
  
  Since both $\ad X$ and $\Ad u$ stabilize the direct sum
  decomposition in Proposition \ref{prop:g2-direct-sum}, and since
  Proposition \ref{prop:g2-regular-elements} shows that both act as a
  single Jordan block on the direct summand $V(\varpi_1)$, the theorem
  follows.
\end{proof}

\subsection{Adjoint partitions}

We  now use the techniques just developed to find
the adjoint partition  for each nilpotent class. We begin with
some lemmas.

Let $L$ be a reductive group with derived group $L' \iso \SL_2(k)$
and suppose that $\dim L = 4$. Let $Z$ be a 1 dimensional central
torus in $L$, and let $W$
be a 4 dimensional simple $L$ module. Fix a maximal torus $S$ of $L'$;
the $S$ weights on $W$ must be $\pm 3, \pm 1$. Since $W$
is simple, $Z$ acts with a single weight $\zeta$ on $W$; suppose
that the integer $\zeta$ is non-0 $\pmod p$.
\begin{lem}
  \label{lem:sl2-orbit}
  If $0 \not = v_1 \in W_{1}$ and $0 \not = v_{-1} \in W_{-1}$, then
  the $L$ orbit of $v =v_1 + v_{-1}$ has dimension 4.
\end{lem}

\begin{proof}
  Let $e_{\pm} \in \Lie(L)_{\pm 2}$ be non-0 weight vectors for $S$.
  Since $W$ is a restricted simple module for $L'$, a theorem of
  Curtis implies that $W$ is simple for $\Lie(L')$; it follows that
  the vectors $e_+v_1 \in W_3$ and $e_-v_{-1} \in W_{-3}$ are
  non-zero.

  Choose $0 \not = h \in \Lie(S)$ and $0 \not = z \in \Lie(Z)$. Then
  $(z,h,e_+,e_-)$ is a basis for $\Lie(L)$.  Let $x = a z + b h + c
  e_+ + d e_-$ with $a,b,c,d \in k$, and suppose that $x.v = 0$.
  
  Since $c e_+ v_1$ is the component of $x.v$ in $W_3$, we deduce that
  $c=0$. Similarly, $d=0$. Thus $x=az + bh$, and $x.v = (a\zeta +
  b)v_1 + (a\zeta -b)v_{-1}$. It follows that $0=a\zeta + b = a\zeta
  -b$; since the characteristic is not 2 and $\zeta$ is non-0 in $k$,
  one has $a=b=0$.
  
  To finish the proof, write $H=\stab_L(v)$. Then $\Lie(H)$ is the
  stabilizer in $\Lie(L)$ of $v$, hence is 0 by the above remarks.
  Thus $\dim H = 0$, so that $\dim L/H = 4$ as desired.
\end{proof}

\begin{rem}
  The proof of the lemma shows that the orbit map $\pi:L \to Lv$ is
  separable.
\end{rem}

\begin{lem}
  \label{lem:g2-distinguished}
  Let $\plie = \plie_{\alpha_1} = \glie_{-\alpha_1} + \lie{b}$, where
  $\blie$ is the Borel subalgebra of $\glie$ corresponding to the
  choice of positive roots determined by $\alpha_1,\alpha_2$. Let
  $\ulie$ be the Lie algebra of the unipotent radical $U$ of the
  parabolic subgroup $P$ with $\Lie(P) = \plie$.  
  Let
  \begin{equation*}
    x_{\alpha_1 + \alpha_2} = [x_{\alpha_1},x_{\alpha_2}] \in  \ulie_{\alpha_1 + \alpha_2}
    \quad \text{and} \quad x_{2\alpha_1 + \alpha_2} = [x_{\alpha_1},x_{\alpha_1 + \alpha_2}] 
    \in \ulie_{2\alpha_1 + \alpha_2}.
  \end{equation*}
  Then 
  \begin{enumerate}
  \item $x = x_{\alpha_1 + \alpha_2} + x_{2\alpha_1 + \alpha_2}$ is a
    representative for the dense (Richardson) orbit of $P$ on $\ulie$.
  \item Put $y_{\gamma_1 + \gamma_2} = [y_{\gamma_1},y_{\gamma_2}],
    y_{\gamma_2 + \gamma_3} = [y_{\gamma_2},y_{\gamma_3}]$, and
    $y_{\gamma_2 + 2\gamma_3} = [y_{\gamma_3},y_{\gamma_2 + \gamma_3}] \in \hlie$.
    Then $d\rho(x) = y_{\gamma_1 + \gamma_2} + y_{\gamma_2 + \gamma_3} + y_{\gamma_2 + 2\gamma_3}$.
  \end{enumerate}

\end{lem}

\begin{proof}
  Let $L = P/U$; thus $L$ is reductive and its derived group $L'$ is simple
  of type $A_1$.  The image of the co-character $\alpha_1^\vee$ is a
  maximal torus of $L'$.  Now, $L$ acts on $\overline\ulie =
  \ulie/[\ulie,\ulie]$ and in fact has a dense orbit on this space.
  Moreover, if $\overline{y} \in \overline{\ulie}$ is a representative
  for the dense $L$ orbit, then any $y \in \ulie$ with image
  $\overline{y}$ lies in the dense (Richardson) $P$-orbit on $\ulie$.
  
  We have $\ulie = \glie_{\alpha_2} + \glie_{\alpha_1 + \alpha_2} +
  \glie_{2\alpha_1 + \alpha_2} + \glie_{3\alpha_1 + \alpha_2} +
  \glie_{3\alpha_1 + 2\alpha_2}$, and $[\ulie,\ulie] =
  \glie_{3\alpha_1 + 2\alpha_2}$. Since $p>3$, the latter assertion
  follows from the commutator formulas in \cite[Theorem 1]{Steinberg};
  these commutator formulas also show that $x_{\alpha_1 + \alpha_2}$
  and $x_{2\alpha_1 + \alpha_2}$ are non-0.
  
  The weights of $\alpha_1^\vee$ on $\overline{\ulie}$ are thus $\pm
  3,\pm 1$; it follows that $\overline{\ulie}$ affords a restricted
  simple 4-dimensional representation for $L'$.  Consider the
  co-character $\gamma = \alpha_1^\vee + 2\alpha_2^\vee$.  Since
  $\langle \alpha_1,\gamma \rangle =0$, the image $Z$ of $\gamma$ is a
  central torus in $L$. Since $\langle \alpha_2,\gamma \rangle = 1$,
  $\gamma$ acts with weight $1$ on $\overline{\ulie}$.  Let
  $\overline{x}$ be the image of $x$ in $\overline{\ulie}$.  Applying
  Lemma \ref{lem:sl2-orbit}, we see that the $L$ orbit through
  $\overline{x}$ has dimension 4, and (1) follows.
  
  For (2) note first that $d\rho(x_{\alpha_1 + \alpha_2}) = y_{\gamma_1
    + \gamma_2} + y_{\gamma_2 + \gamma_3}$.  Since 
  $$d\rho(x_{2\alpha_1
    + \alpha_2}) = [y_{\gamma_1},y_{\gamma_2 + \gamma_3}] +
  [y_{\gamma_3},y_{\gamma_1 + \gamma_2}] + [y_{\gamma_3},y_{\gamma_2 +
    \gamma_3}]$$
  an application of the Jacobi identity shows that $d\rho(x_{2\alpha_1
    + \alpha_2}) = [y_{\gamma_3},y_{\gamma_2 +
    \gamma_3}]$ and (2) follows.
\end{proof}

\begin{lem}
  \label{lem:rep-ring-calcs}
  Let $\formal$ be the additive formal group law, and let
  $R=R_\formal$ be the its representation ring (see
  \S\ref{sub:formal-groups}). For $p > 3$,  we have the following identities in $R$:
  \begin{enumerate}
  \item $\bwedge^2(2J_3 + J_1) = J_5 + 5J_3 + J_1$.
  \item $\bwedge^2(J_3 + 2J_2) = 2J_4 + 2J_3 + 2J_2 + J_1$.
  \item $\bwedge^2(2J_2 + 3J_1) = J_3 + 6J_2 + 6J_1$.
  \item If $p=5$ or $p \ge 11$, $\bwedge^2(J_7) = J_{11} + J_3$.
    If $p=7$, $\bwedge^2(J_7) = 2J_7$.
  \end{enumerate}
\end{lem}

\begin{proof}
  The first three assertions are an immediate consequence of the
  following facts which the reader may easily check: $\bwedge^2 J_2 = J_1$,
  $\bwedge^2 J_3 = J_3$,  $J_3 \tensor J_3 = J_5 + J_3 + J_1$,  $J_3
  \tensor J_2 = J_4 + J_2$, and  $J_2 \tensor J_2 = J_3 + J_1$.

  For $p \ge 7$, (4) follows by considering the tilting
  $\SL_2(k)$-module $\bwedge^2 L(6)$, and applying Proposition
  \ref{prop:tilting-partitions}. When $p=5$ one checks (4) by hand (or
  by computer).
\end{proof}

\begin{theorem}
  \label{theorem:g2-partitions}
  The partitions of the non-0 nilpotent classes of $\glie$ on
  $V(\varpi_1)$ and on the adjoint representation are as follows (recall that $p>3$):
\begin{equation}
  \label{eq:g2-partitions}
  \begin{array}[c]{c|lll}
    \text{nilpotent orbit} & \text{partition on $V(\varpi_1)$} & \text{adjoint partition} \\
    \hline \hline
    A_1 & (2^2,1^3) & (3,2^4,1^3) \\
    \widetilde{A_1} & (3,2^2) & (4^2,3,1^3) \\
    G_2(a_1) & (3^2,1) & (5,3^3) \\
    G_2 & (7) &  (11,3) & \text{if $p \ge 11$ or $p=5$} \\
    & & (7^2) & \text{if $p=7$.} \\
  \end{array}
\end{equation}

\end{theorem}

\begin{proof}
  The adjoint partitions may be obtained from the partitions on $V(\varpi_1)$ by
  applying Lemma \ref{lem:rep-ring-calcs} and Proposition
  \ref{prop:g2-direct-sum}.
  
  To obtained the partitions on $V=V(\varpi_1)$ note the following.
  $x_{\alpha_2}$ is a representative for the class $A_1$, and
  $x_{\alpha_1}$ is a representative for the class $\widetilde{A_1}$;
  one now deduces the partitions from \eqref{eq:x-a-1} and
  \eqref{eq:x-a-2}.  The partition of the regular class $G_2$ is
  obtained from Proposition \ref{prop:g2-regular-elements}.
  
  Lemma \ref{lem:g2-distinguished} gives a representative $x$ for the
  class $G_2(a_1)$; moreover, it shows that $d\rho(x) = a + b + c$
  with $0 \not = a \in \hlie_{\gamma_1+\gamma_2}$, $0 \not = b \in
  \hlie_{\gamma_2+\gamma_3}$ and $0 \not = c \in
  \hlie_{\gamma_2+2\gamma_3}$.  If $a,b,c$ are \emph{any} non-0 elements of
  the indicated weight spaces of $\hlie$, a direct calculation shows that
  the partition of $a+b+c$ on $V$ is $(3,3,1)$.
\end{proof}

\begin{rem}
  According to Theorem \ref{theorem:g2-adjoint-partitions} the
  partitions in \eqref{eq:g2-partitions} are
  also valid for the unipotent classes. Lawther \cite{Lawther:Jblock}
  has computed adjoint Jordan blocks for unipotent elements in
  exceptional groups, and our results agree with his in this case.
  
  The descriptions in \S\ref{sub:simple-roots-g2} imply that the
  partitions of the unipotent classes on $V(\varpi_1)$ are given by
  \eqref{eq:g2-partitions}; this again agrees with Lawther's
  calculations.
\end{rem}

\section{Adjoint Jordan blocks in characteristic 0}
\label{sec:char0}

In this section, we work over an algebraically closed field of
characteristic 0; in order to emphasize the characteristic, we will
call this field $F$ rather than $k$.  $G$ will be a simple algebraic
group over $F$, and $\glie = \Lie(G)$ will be its Lie algebra.

Fix a distinguished nilpotent element $X \in \glie$, and choose an
$\sllie_2$-triple in $\glie$ containing $X$, and let $\phi:\SL_2(F)
\to G$ be the homomorphism such that the image of $d\phi$ is the span
of this $\sllie_2$-triple.  We may suppose that the cocharacter
$\nu(t) = \phi(\operatorname{diag}(t,t^{-1})):F^\times \to G$
satisfies $\Ad \nu(t)X = t^2X$ ($t \in F^\times$).

Let $\glie(i) = \glie(i;\nu)$ denote the $i$-th weight space for $\Ad
\nu$.  Since $X$ is distinguished, \cite[5.7.6]{Carter1} shows that
$\glie(2i+1) = 0$ for $i \in \Z$ .  

Let $n \ge 0$ be the largest integer for which $\glie(2n) \not = 0$.
Fix $\zeta \in F^\times$ a root of unity of order $n+1$, and
let $s = \nu(\zeta^{1/2}) \in G$ (choose either square root of $\zeta$
in $F$). The eigenvalues of $\Ad(s)$ on $\glie$ are
integeral powers of $\zeta$; the $\zeta^i$ eigenspace
of $\Ad(s)$ is
\begin{equation*}
  V_i = \glie(2i) \oplus \glie(2i-2n-2) \quad \text{for} \ 0 \le i \le n.
\end{equation*}

Let $M \in \glie(-2n)$, and put $S = X + M \in V_1$.
\begin{lem}[Springer]  \cite[Lemmas 9.3, 9.5, and 9.6]{springer-regelts}
  \label{lem:springer}
  \begin{enumerate}
  \item $\Ad(s)$ stabilizes $C_\glie(S)$ and has no non-0 fixed points
    in $C_\glie(S)$.
 % \item $S$ is nilpotent if and only if $M = 0$.
%  \item $C_G(S)^o$ is a solvable group.
  \item $V_1$ contains a regular semisimple element if and only if 
    $S$ is regular semisimple for some choice of $M \in \glie(-2n)$.
  \item If $\dim \glie(4) = \dim \glie(2)-1$, then $V_1$ contains a
    regular semisimple element.
  \end{enumerate}
\end{lem}

Let $\nlie = \bigoplus_{i > 0} \glie(i)$ and $\nlie^- = \bigoplus_{i <
  0} \glie(i)$.  The following lemma may be found in \cite[Lemma
6.4A]{kostant-betti} for the case where $X$ is regular; we have
essentially copied Kostant's argument.
\begin{lem}
  \label{lem:rho-injective}
  Let $\rho:\glie \to \nlie$ be the projection with kernel $\nlie^-
  \oplus \glie(0)$.  Then $\rho \circ \Ad(s) = \Ad(s) \circ \rho$, 
  the restriction $\rho_{\mid \clie}$ to $\clie = \clie_\glie(S)$ takes values in
  $\clie_\glie(X)$, and $\rho_{\mid \clie}$ is injective.
\end{lem}

\begin{proof}
  The $\Ad(s)$-equivariance of $\rho$ is immediate.
  Let $Y \in \clie$ and write
  \begin{equation*}
    Y = U + A + V \in \nlie^- \oplus \glie(0) \oplus \nlie.
  \end{equation*} 
  We claim  that $V = \rho(Y) \in
  \clie_\glie(X)$.  Indeed, note that $$[M,Y] = [M,A] + [M,V] \in
  \nlie^- + \glie(0),$$ so that $\rho([M,Y]) = 0$.  Since $Y \in \clie$, one has
  $0=\rho([S,Y]) = \rho([X,Y] + [M,Y]) = \rho([X,Y])$. Since $[X,U]
  \in \nlie^- + \glie(0)$, we have 
  $$0=\rho([X,Y]) = \rho([X, A]) + \rho([X,V]) \in \glie(2) \oplus
  \bigoplus_{i \ge 4} \glie(i). $$
  It follows that $[X,A] =
  \rho([X,A]) = 0$ and $[X,V] = \rho([X,V]) = 0$; this proves that
 $V \in \clie_\glie(X)$.
  
 Since $X$ is distinguished, $\ad(X):\glie(0) \to \glie(2)$ is
 bijective \cite[5.7.4]{Carter1}; it follows that $\clie_\glie(X)
 \subset \nlie$ (view $\glie$ as a representation for the
 $\lie{sl}_2(F)$-triple to see that $\clie_\glie(X) \subset \glie(0) +
 \nlie$).  Since $[X,A]=0$, we have $A=0$. Thus we see that $Y \in
 \clie$ has the form $Y = U + V \in \nlie^- \oplus \nlie^+$.
 
 Now suppose that $Y \in \ker \rho \cap \clie$; by the previous
 remarks we have $Y \in \nlie^-$.  Considering homogeneous components
 of $Y$, one sees that $[M,Y] = 0$.  It follows that $0=[S,Y] = [M,Y]
 + [X,Y] = [X,Y]$ i.e. $Y \in \clie_\glie(X)$.  Since $\clie_\glie(X)
 \cap \nlie^- = 0$, we have $Y = 0$ and the lemma is proved.
\end{proof}

Let $\hlie = \Lie(T)$ be the Lie algebra of a maximal torus of $G$
with Weyl group $W = N_G(T)/T$. Then $W$ acts by graded algebra
automorphisms on the coordinate ring $F[\hlie]$.  Since $W$ is a
reflection group, a theorem of Chevalley says that the algebra of
invariants $F[\hlie]^W$ is generated by $r = \dim \hlie$ algebraically
independent homogeneous polynomials $\phi_1,\dots,\phi_r$ whose
degrees are uniquely determined, up to order.  The exponents of $W$
are the numbers $e_i = \deg \phi_i - 1$, $1 \le i \le r$.

The following theorem generalizes the result proved by Kostant
\cite{kostant-betti} in the case where $X$ is regular nilpotent.

\begin{theorem}
  \label{theorem:char-0-jblocks}
  Suppose that $V_1$ contains a regular semisimple element of $\glie$.
  Let $e_1,\dots,e_r$ be the exponents of the Weyl group of $G$, and
  for $1 \le i \le r$, let $1 \le f_i \le n$ be the unique quantity
  satisfying $f_i \congruent -e_i \pmod n$.  Then $\ad X$ has a Jordan
  block of size $2f_i + 1$ for each $1 \le i \le r$.
\end{theorem}

\begin{proof}
  Since $V_1$ contains a regular semisimple element, Lemma
  \ref{lem:springer}(2) shows that we may choose $M \in \glie(-2n)$
  such that $S = M + X$ is regular semisimple. Then the centralizer
  $\hlie' = \clie_\glie(S)$ is the Lie algebra of a maximal torus $T'$
  of $G$.
  
  Now part (1) of Lemma \ref{lem:springer} shows that
  $s=\nu(\zeta^{1/2})$ determines an element $w$ of the Weyl group
  $N_G(T')/T'$. Moreover, $w$ has a regular eigenvector in
  $\hlie'$, namely $S$; it is thus a \emph{regular}
  Weyl group element; see \cite[\S4]{springer-regelts} for these
  notions. In particular, \cite[Theorem 4.2]{springer-regelts} shows
  that the eigenvalues of $w$ (i.e. of $\Ad(s)$) on $\hlie'$ are $\zeta^{-e_i} =
  \zeta^{f_i}$ for $1 \le i \le r$.
  
  Let $\lie{s}$ be the $\lie{sl}_2(X)$-triple containing $X$; as
  $\lie{s}$-module, we may write $$\glie = L(2\lambda_1) \oplus
  L(2\lambda_2) \oplus \cdots \oplus L(2\lambda_t).$$
  Here $L(\mu)$ is
  the simple $\lie{s}$-module of dimension $\mu+1$; since
  $S$ is distinguished, any simple
  $\lie{s}$-summand of $\glie$ has even highest weight.
  
  We may chose $v_i \in \glie(2\lambda_i)$ for $1 \le i \le t$ such
  that $v_1,\dots,v_t$ is a basis for $\clie_\glie(X).$ It follows
  that the eigenvalues of $\Ad(s)$ on $\clie_\glie(X)$ are
  $\zeta^{\lambda_j}$ for $1 \le j \le t$.
  
  By Lemma \ref{lem:rho-injective} the map $\rho:\clie_\glie(S) \to
  \clie_\glie(X)$ is injective and $\Ad(s)$-equivariant. Re-ordering
  the $\lambda_i$ if necessary, we may therefore suppose that
  $\zeta^{\lambda_j} = \zeta^{f_j}$ for $1 \le j \le r$.  For each $j$
  we have $1 \le \lambda_j \le n$ and $1 \le f_j \le n$, whence
  $\lambda_j = f_j$ and the theorem follows.
\end{proof}

\begin{rem}
  The hypothesis that $V_1$ contains a regular semisimple element is
  satisfied if the condition in (3) of Lemma \ref{lem:springer} holds.
  In \cite{springer-regelts}, Springer lists the nilpotent classes for
  which (3) holds (though there is class for type $G_2$ in that list
  that doesn't belong). 
  
  Springer's results yield a map from nilpotent orbits for which $V_1$
  contains a regular semisimple element to conjugacy classes in the
  Weyl group. Kazhdan and Lusztig \cite{KL-affine-flag} have defined a
  map from all nilpotent orbits to conjugacy classes in the Weyl
  group; their map agrees with Springer's when his is defined.
  
  For a nilpotent $X \in \glie$, let $\sigma(X)$ denote the
  corresponding conjugacy class in $W$ under the Kazhdan-Lusztig map.
  In the case where $X$ is distinguished and $\sigma(X)$ consists of
  regular elements of $W$ in Springer's sense, it appears from
  empirical observation that Theorem \ref{theorem:char-0-jblocks} is
  still valid, i.e.  one obtains $r$ of the Jordan block sizes of
  $\ad(X)$ from the eigenvalues of $\sigma(X)$. I have not so far been
  able to give a proof which explains this phenomenon.
\end{rem}

%\bibliography{MathBib,extra-bibs}
\providecommand{\bysame}{\leavevmode\hbox to3em{\hrulefill}\thinspace}

\end{document}